\documentclass[12pt]{article}
\usepackage{amssymb, amsmath, amsfonts}
\def\nodo#1{}

\def\max{\operatorname{max}}
\font\tenDDl=msbm10  
\font\sevenDDl=msbm7 
\font\fiveDDl=msbm5 
        \newfam\DDlfam  
                \textfont\DDlfam=\tenDDl
\scriptfont\DDlfam=\sevenDDl
        \scriptscriptfont\DDlfam=\fiveDDl
\begin{document}
\begin{center}
{\bf \Large Every quantum minor generates an Ore set}
\vskip .2in
{\sc Zoran \v{S}koda} {\tt zskoda@irb.hr}

{\bf Abstract.} {\it The subset multiplicatively generated by
any given set of quantum minors and the unit element
in the quantum matrix bialgebra 
satisfies the left and right Ore conditions.
}
\end{center}
Quantum matrix groups~(\cite{ManinQGNG,ParshallWang, RTF}) 
have remarkable algebraic properties and
connections to several branches of mathematics and mathematical physics.
A viewpoint of the noncommutative geometry
may elucidate some of their properties. 
In the formalism in which a quantum group is
described by a matrix bialgebra
${\cal G}$, e.g.  ${\cal SL}_q(n)$,
the shifts of the main Bruhat cells are expected also to have
quantum analogues which are localizations of ${\cal G}$. 
The geometry is richer and more akin to the classical case, 
if these localizations have good flatness properties. 
Ore localization is the most well-understood kind of
localizations. Ore localizations are biflat 
(in terminology of~\cite{Ros:NcSch}), 
and appear often in ``quantum'' situation, that is, when the
noncommutative algebra is just a deformation of a commutative
algebra.   
Thus one expects to realize the quantum main Bruhat cell and
its Weyl group shifts as the (noncommutative spectra of) 
the localized algebras of the form ${\cal G}[S^{-1}_w]$
where $S_w$ are certain Ore subsets in ${\cal G}$, 
depending on the element $w$ in the Weyl group $W$. 
Furher support of this conjectural picture is a result of {\sc A.~Joseph}
~\cite{Joseph:fflat}, that there is a natural family 
of Ore sets $\bar{S}_w$, $w \in W$, 
in the graded algebra ${\cal R}$ representing the quantum analogue of the
basic affine space $G/U$, such that in the commutative case
the spectra of the localizations 
${\cal R}[\bar{S}^{-1}_w]$ are exactly the images
of the Bruhat cells in the basic affine space. 

However, the Ore property for the natural candidates for ${\cal S}_w$
has not been proved so far. Trying to answer
the question of {\sc V.~Drinfeld} to find the quantum analogue of
the {\sc Beilinson-Bernstein} localization theorem, 
{\sc Y. Soibelman} has shown the satisfactory localization 
picture for ${\cal SL}_q(2)$ (unpublished), 
and came to the conclusion that for ${\cal SL}_q(n)$, $n>2$, 
already the Ore property of $S_w$ is far from obvious, if not even wrong.
Quantum Beilinson-Bernstein theorem has been further studied by
{\sc Lunts} and {\sc Rosenberg}~\cite{LuntsRos1}, 
{\sc Tanisaki} \cite{Tanisaki:BeilBern} and others,
but in different approaches.

In his thesis, the present author has 
developed a direct localization approach
~(\cite{Skoda:ban, Skoda:qGauss}) to the 
construction of the coset spaces 
of the quantum linear groups and the locally trivial quantum 
principal fibrations deforming the classical fibrations
of the type $G \to G/B$ and having Hopf algebras as
the replacements of the structure groups. 
Apart from the sketch in~\cite{Skoda:ban}, 
the main part of that work has been still unpublished
(however, a nontrivial application to quantum group 
coherent states and appropriate measure is exhibited in~\cite{Skoda:coh}). 
The present paper fills a part of this gap in view of the
observation that the sets $S_w$ in ${\cal SL}_q(n)$ 
are sets multiplicatively generated by a specific set of quantum minors 
attached to the permutation matrix $w$, 
namely the set of all principal (=lower right corner)
quantum minors of the row-permuted matrix of generators 
$T = (t^{w^{-1}(i)}_j)$. 
Note that this solves the problem
weather the quantum Bruhat cells are realizable as Ore localizations only 
for quantum matrix groups of type A (in the sense of Lie theory). 

A subset $S$ in a ring $R$ is {\bf multiplicative}
if $1 \in S$ and $s,s'\in S$ implies $ss' \in S$. 
A subset $S\in R$ satisfies the {\bf left Ore condition} if
$\forall r\in S, \forall r\in R$, $\exists r'\in R$ $\exists s'\in S$
with $s' r = r's$. If $R$ has no zero divisors, then 
a multiplicative subset $S\subset R$ satisfying the left Ore condition
is called a {\bf left Ore set}. 
The right Ore condition is the left Ore condition in the
ring with opposite multiplication $R^{\rm op}$. 
A subset $S \in R$ is Ore 
if it is left Ore in $R$ and $R^{\rm op}$ simultaneously.
For more details on Ore sets see~\cite{Skoda:nloc, vOyst:assalg}.

Let $I$, $J$ be some linearly ordered sets of equal cardinality 
$|I| = |J| = n >0$, where the elements of $K$ and $L$ 
are called labels. Given a field ${\bf k}$ of characteristic zero,
and number $0\neq q \in {\bf k}$,
the underlying ${\bf k}$-algebra of the
matrix bialgebra ${\cal M}_q(I,J)$ is the free ${\bf k}$-algebra
on generators $t^i_j$ where $i \in I$, $j \in J$, 
modulo the relations
\[\begin{array}{l}
 t^\alpha_\gamma t^\alpha_\delta = q t^\alpha_\delta t^\alpha_\gamma,
\,\,\,\,\,\,\,\,\,\,\alpha = \beta, \,\,\gamma < \delta \,\mbox{ (same row) }
\\
 t^\alpha_\gamma t^\beta_\gamma = q t^\beta_\gamma t^\alpha_\gamma,
\,\,\,\,\,\,\,\,\,\,\alpha < \beta, \,\,\gamma = \delta 
\,\mbox{ (same column) }
\\ 
t^\alpha_\gamma t^\beta_\delta - t^\beta_\delta t^\alpha_\gamma= 
 \left\lbrace
\begin{array}{l}
(q - q^{-1}) t^\beta_\gamma t^\alpha_\delta ,\,\,\,\,\,
 \alpha < \beta,\,\,\gamma < \delta\\
0 ,\,\,\,\,\,\,\,\,\,\,\,\,\,\,\,
 \alpha < \beta,\,\, \gamma > \delta\end{array}\right.
\end{array}\]
The isomorphism class of algebra ${\cal M}_q(K,L)$
depends only on the cardinality $n$ of $K$ and $L$, and we denote
a representative of that class ${\cal M}_q(n)$, e.g. 
where $K = L = {\bf n} = \{1,\ldots, n\}$.

Let $K = (k_1,\ldots, k_m)\subset I$, $L = (l_1,\ldots, l_m)\subset J$ 
be subsets of equal cardinality $m<n$
where the labels in the round brackets are ordered
according to the subset order from $I$ and $J$. 
Then $T^K_L$ will be the submatrix of
$T = (t^i_j)^{i \in I}_{j \in J}$ whose rows
have labels in $K$ and whose columns have labels in $L$.
Denote $D^K_L = {\rm det}_q T^K_L = \sum_{\sigma \in \Sigma(m)}
(-q)^{l(\sigma)} t^{k_1}_{l_{\sigma(1)}} \cdots t^{k_m}_{l_{\sigma(m)}}$
where $l(\sigma)$ denotes the number of inversions in the permutation
$\sigma$ on the set of $m$ labels. Any element of this form 
is called a {\bf quantum minor}.  
$D^K_L$ is central in ${\cal M}_q(K,L)$, but
not in the whole ${\cal M}_q(I,J)$, unless $K = I$ and $L = J$.

{\bf Lemma 1.} {\it Let $S\subset R$ be a multiplicative subset
and $E = E(S)$ be the set of all $e \in R$ 
which satisfy the following 'partial' left Ore condition
\begin{equation}\label{eq:partialOre} 
\forall s \in S, \, \exists s' \in S, \exists r' \in R \,\,\,\,\,
r' s = s' e 
\end{equation}
Suppose also that $S \subset E$. Then (i) $E$ is a subalgebra of $R$.

(ii) If $r \in R$ satisfies the following $E$-relative partial 
left Ore condition
\begin{equation}\label{eq:qmin-ore-proof-1}
 \forall s \in S, \, \exists s' \in S, \exists r' \in R \,\,\,\,\,\,
r's - s' r \in E \end{equation}
then $r \in E$.

(iii) Let $S_0 \subset S$ multiplicatively generates $S$.
Then for left Ore condition to hold, that is $E(S) = R$, 
it is enough to check~\eqref{eq:partialOre} 
for all $e\in R$ but only $s\in S_0$ (still with $s' \in S$).
}

{\it Proof.} The reasoning method here is pretty standard 
(cf.~\cite{Skoda:nloc}, Ch. 6).

(i) Given $e_1,e_2 \in E$, we prove that $e_1 e_2 \in E$ as follows. 
By the assumption, given any $s\in S$, 
$\exists r'\in R$, $\exists s'\in S$, with  $r' s = s' e_1$;
and $\exists r''\in R, \exists s''\in S$ such that $r'' s' = s'' e_2$.
Then $s'' (e_1 e_2) = r'' s' e_2 = (r'' r') s$. 

To prove that $e_1 + e_2 \in E$ we reason as follows. 
Given any $s \in S$, we first find $s_1,s_2 \in S, r_1,r_2\in R$
such that $s_1 e_1 = r_1 s$ and $s_2 e_2 = r_2 s$. 
By $S \subset E$, also $\exists s_* \in S$, $r_* \in R$ such that
$s_* s_1 = r_* s_2$. Then
$$\begin{array}{lcl}
(s_* s_1) (e_1 + e_2) &=& s_* r_1 s + s_* s_1 e_2\\ 
                      &=& s_* r_1 s + r_* s_2 e_2 \\
                      &=& s_* r_1 s + r_* r_2 s \\  
                      &=& (s_* r_1 + r_* r_2)s,\,\,\,\,\,\,\,\,\mbox{with }
s_* s_1\in S,\mbox{ as required.}
\end{array}$$ 

(ii) Let $r' s - s' r = e \in E$. By $e \in E$ we can find 
$r'', s''$ such that $r''s = s'' e$. 
Then $(s'' r' - r'')s = s''(r's - e) = s'' s' r$ 
with $s'' s' \in S$, as required. 

(iii) Suppose $s_1,s_2 \in S_0$. By assumption, 
$r_2 s_2 = s_2' e$ and $s_* r_2 = r_* s_1$ 
for some $s_2', s_* \in S$ and some $r_2,r_* \in R$. 
Then $r_* (s_1 s_2) = s_* r_2 s_2 = (s_* s_2') e$ with
$s^* s_2' \in S$. Hence, the equation~\eqref{eq:partialOre}
holds for all $s\in S'_0$ where $S'_0$ consists of all products
of pairs of elements in $S$. Continuing by induction on
the length $n$ of a product $s = s_1\ldots s_n \in S$ we conclude
that~\eqref{eq:partialOre} holds for all $s\in S$.  

{\bf Lemma 2.} {\it 
Let $|K|=|L| < n$.  Consider the subalgebra $E_0 = E_0(K,L)$
in ${\mathcal M}(n)$ 
generated by all $t^k_l$ where either $k \in K$ or $l\in L$ (or both). 
Let $t^{k'}_{l'} \notin E_0$. Then for every $e \in E_0$,
$t^{k'}_{l'} e - e t^{k'}_{l'} \in E_0$.
}

{\it Proof.} This is a linear statement, hence it suffices to
prove it for words $e = t^{k_1}_{l_1}\cdots t^{k_r}_{l_r}$ of degree $r$
in generators of $E_0$. We show this by induction on $r$. 
For $r = 0,1$, this is obvious. 
Let $e = e_{r-1} t^{k_r}_{l_r}$ where $e_{r-1}$ is of the degree $r-1$. 
The commutator $[e_{r-1}, t^{k'}_{l'}]\in E_0$ by the induction hypothesis,
and $[t^{k_r}_{l_r}, t^{k'}_{l'}]$ is either zero or proportional to 
$t^{k'}_{l_r} t^{k_r}_{l'}$. Since $E_0$ is a subalgebra, these two facts
and the identity $[e_{r-1} t^{k_r}_{l_r}, t^{k'}_{l'}] = 
[e_{r-1}, t^{k'}_{l'}] t^{k_r}_{l_r} + e_{r-1} [t^{k_r}_{l_r}, t^{k'}_{l'}]$
imply that $[e,t^{k'}_{l'}]\in E_0$.\,\,\,Q.E.D.

{\bf Notation.} If $L = (l_1,\ldots, l_r)$, then 
$L(l_k\to l) := (l_1,\ldots, l_{k-1}, l, l_{k+1},\ldots,l_r)$.

{\bf Lemma 3.} (special cases of commutation relations for quantum minors)

(i) ~~\cite{ParshallWang}
If $k \in K$ and $l \in L$ then $t^k_l D^K_L = D^K_L t^k_l$.  

(ii) ~\cite{ParshallWang}
If either $k > \max K$ and $l \in L$, or $l > \max L$ and $k \in K$, 
then $t^k_l D^K_L = q^{-1} D^K_L t^k_l$. 
If either $k < \max K$ and $l \in L$, or $l < \max L$ and $k \in K$, 
then $t^k_l D^K_L = q D^K_L t^k_l$. 
Also (say, by Muir's law~\cite{KrobLec}),
if $L$ and $L'$ differ by an interchange of a single label, 
then $D^K_L D^K_{L'} = q^{\pm 1} D^K_{L'} D^K_L$. 

(iii) Let $L = (l_1,\ldots, l_s)$. Suppose $l_1 < l < l_i$, for $i>1$. Then
\begin{equation}\label{eq:qmin-ore-proof-2}
D^K_L t^k_l - q^{-1} t^k_l D^K_L = 
(1-q^{-2}) D^{K}_{L(l\mapsto l_1)} t^{k}_{l_1}
\end{equation}

(iv) Let $L = (l_1,\ldots, l_s)$ be ascendingly ordered and
$l_r < l < l_{r+1}$ for some $r$. Then (iii) generalizes to
\begin{equation}\label{eq:r-gap-formula-in}
D^K_L t^k_l - q^{-1}t^k_l D^K_L
 = q^{-1}(q-q^{-1}) 
\sum_{u=1}^r (-q)^{u-r} t^{k_m}_{l_u} D^K_{L(l_u \rightarrow l)}
\end{equation}
Parts (iii) and (iv) are now widely known among people working 
on ``identities between quantum minors'', 
and were many times independently rediscovered 
by many people (e.g. by the author around 1998). Both identities are
much easier to prove when $k = \max{K}$ (using for example Laplace
identities~\cite{Hay:qmultilinear, KrobLec, ParshallWang}, 
and simple calculational arguments, and, for (iv), 
induction). The general case, follows from $k = \max{K}$ case 
applying the included-row exchange principle~(\cite{Skoda:inclrow}).

{\bf Fact.} (e.g.~\cite{ParshallWang}) 
{\it ${\cal M}_q(n)$ has no zero divisors.}

{\bf Theorem.} {\it Each quantum minor $D^K_L = {\rm det_q}\, T^K_L$ 
(where $K$ and $L$ are the subsets of $\{1,\ldots,n\}$ of the same cardinality)
and $1$ multiplicatively generate an Ore set $S = S^K_L$
in $R = {\cal M}_q(n)$.
}

{\it Proof.} We prove the right Ore condition. 
The left Ore condition then follows
by using the algebra automorphism $t^i_j \rightarrow t^j_i$.

According to the Lemma 1 (i),  
it is sufficient to prove that all the generators $t^k_l$ 
are in $E = E(S^K_L)$, and $S \subset E$. The latter is clear
as $S$ is multiplicatively generated by only one element.
Furthermore, by Lemma 1 (iii) we will be checking 
the condition~\eqref{eq:partialOre} for $s = D^K_L$ only,
and not for the higher powers of $D^K_L$.
We first prove for $t^k_l \in E_0(K,L)\subset R$, (in notation of Lemma 2). 

If $k \in K$ and $l \in L$ this is obvious because then
$t^k_l$ and $D^K_L$ commute. If $k \in K$ and $l$ is smaller than
the minimum of $L$ or larger then the maximum of $L$
then $t^k_l D^K_L = q^{\pm 1} D^K_L t^k_l$ what again proves the condition.
The same way conclude that in the case that $l \in L$ and $k$ is
smaller than the minimum of $K$ or larger than the maximum of $K$.

Now consider the case when $k \in K$ and $l_r < l < l_{r+1}$ for some
$r$. If $r = 1$ then the formula~\eqref{eq:qmin-ore-proof-2} applies.
Notice that the RHS is in $E$, because $t^k_{l_1}$ is
already proven to be in $E$, and $D^K_L$ and $D^K_{L'}$
differ only in one column, hence $q$-commute according to
Lemma 3 (ii) and $E$ is a subalgebra.
Hence the formula ~\eqref{eq:qmin-ore-proof-2} is easily recognized
to be of type ~\eqref{eq:qmin-ore-proof-1} and therefore $t^k_l \in E$.

Now we use induction on $r$ to prove the same for all $r$.
So suppose we have proved for $r-1$.
Then, the formula~\eqref{eq:r-gap-formula-in} applies.
Using the induction hypothesis, one easily sees 
that the factors involved in each summand
on the right-hand side of~\eqref{eq:r-gap-formula-in} are in $E$:
Namely,  $D^K_{L(l_u \rightarrow l)}$ has at most $(r-1)$ columns with indices
less than $l$, and $t^{k_m}_{l_u}$ is one of the entries in $T^K_L$,
hence it commutes with $D^K_L$.

Thus we have proved the claim for all $t^k_l \in E_0 = E_0(K,L)$,
and therefore $E_0 \subset E$.  
We are left with the case of $t^k_l\notin E_0$, i.e. where both
$k \notin K$ and $l \notin L$. Surely $D^K_L \in E_0$, hence by
Lemma 2, $t^k_l D^K_L - D^K_L t^k_l \in E_0\subset E$. 
By Lemma 1 (ii) $t^k_l \in E$. Hence by Lemma 1 (i) we are done.

{\bf Corollary.} {\it (0) The theorem holds for 
${\mathcal SL}_q(n)$ and ${\mathcal GL}_q(n)$.

(i) Every set of quantum minors
multiplicatively generates an Ore set in $R= {\mathcal M}_q(n),
{\mathcal SL}_q(n)$ or ${\mathcal GL}_q(n)$.

(ii) The product of finitely many $q^r$-commuting (various $r$-s)
quantum minors is a multiplicative generator of an Ore set in $R$.

(iii) The product of all principal quantum minors
(lower right corner) of sizes $1 \times 1$ to $(n-1)\times (n-1)$ 
of any row and column-permuted matrix
$G = T^\sigma_\tau$ of $T$ is a multiplicative generator 
of an Ore set in $R$.

(iv)  The theorem holds for the multiparametric quantum groups 
obtained by twisting~
(\cite{Artin:multi, Demidov:rev, Manin:multi, Skoda:multipar}).
}

To observe (0), notice that ${\mathcal SL}_q(n))$ has no zero divisors and 
it is a quotient of  ${\mathcal M}_q(n)$. The projection map
sends all quantum minors to nonzero elements, so the Ore
conditions are fullfilled by using the projection map directly. 
${\mathcal GL}_q(n))$ is a {\it central} localization of ${\mathcal M}_q(n)$,
hence compatible with any Ore set in ${\mathcal M}_q(n)$: any central
localization sends Ore sets to Ore sets. 

(i-iii) easily follow from the theorem by general principles on 
recognizing Ore sets~(\cite{Skoda:nloc}, Ch.~6). 

(iv) It is not difficult to observe, and it was shown in
detail in~\cite{Skoda:multipar}, the twisting of~\cite{Artin:multi}
changes the relations (in Lemma 3) among the quantum minors by nonzero 
factors in each of the summands, and the rest of the reasoning follows
the same way, up to numerical factors, which can always be absorbed in
$r'$ in Ore condition. 

In fact, one can directly use the theorem in its 1-parameter version, 
and an isomorphism from the multiparametric quantum bialgebra into the
$(0,0)$-bidegree component of certain tensor product 
$S_l \otimes {\cal M}_q(n) \otimes S_r$. 
This isomorphism appropriately transfers the
quantum minors and Ore conditions involving certain class of 
homogeneous elements in general~(\cite{Skoda:multipar}). 
\vskip .1in

{\footnotesize
{\bf Acknowledgements.} The main theorem has been proved by the author in
February 2000, and presented among the basic lemmas in his Wisconsin
thesis (defended January 2002), and anounced in~\cite{Skoda:ban}, Theorem 8. 
I thank Prof. {\sc J. W. Robbin} for his mathematical mentorship,
and Prof. {\sc Y.~Soibelman} for a discussion and his
assesment of the importance of the result in Madison, 2001.   
The exposition has been adapted to the present form 
at Max Planck Institute for Mathematics, Bonn,
and the final version at IRB, Zagreb. 
}

\end{document}